\documentclass[12pt,twoside,draft]{article}

\usepackage{amsmath,amssymb}

\numberwithin{equation}{section}

\begin{document}

\centerline{}

\centerline{}

\centerline {\textsc {\Large A generalization of}}
\centerline {\textsc {\Large random self-decomposability}}

\centerline{}

\centerline{\textsf {\large Satheesh S}}

\centerline{NEELOLPALAM, S. N. Park Road}

\centerline{Thrissur-680 004, India.}

\centerline{e-mail: \textit{ssatheesh1963@yahoo.co.in}}

\centerline{}

\centerline{\textsf {\large Sandhya E}}

\centerline{Department of Statistics, Prajyoti Niketan College}

\centerline{Thrissur-680 301, India.}

\centerline{e-mail: \textit{esandhya@hotmail.com}}

\centerline{}

{\bf Abstract.} The notion of random self-decomposability is generalized here. Its relation to self-decomposability, Harris infinite divisibility and its connection with a stationary first order generalized autoregressive model are presented. The notion is then extended to $\mathbf{Z_+}$-valued distributions.



\section{\large Introduction}\label{sec1}

The role of self-decomposable (SD) distributions in first order autoregressive (AR(1)) models of the form
\begin{equation}
X_n = cX_{n-1}+\epsilon_n,
\end{equation}
\noindent described by random variables (\textit{r.v.}s) $\{X_n, n\in Z\}$, innovations (\textit{i.i.d. r.v.}s) $\{\epsilon_n\}$ and $c\in (0,1)$ such that for each $n$, $\epsilon_n$ is independent of $X_{n-1}$, has been discussed by many authors, \textit{see e.g.} Bouzar and Satheesh (2008) and the references therein. Recently Kozubowski and Podg\'orski (2010) has introduced the notion of random self-decomposability of distributions on the reals motivated by stationary solutions to the AR(1) model

\vspace{2mm}
   
\begin{equation}
X_n = \begin{cases}
\epsilon_n, \text{ with probability $p$}, \\
cX_{n-1}+\epsilon_n, \text{ with probability $(1-p)$,}
\end{cases}
\end{equation}
\noindent described by \textit{r.v.}s $\{X_n, n\in Z\}$, innovations $\{\epsilon_n\}$ and $c\in [0,1]$ such that for each $n$, $\epsilon_n$ is independent of $X_{n-1}$.

\vspace{2mm}

\noindent \textbf {Definition 1.1} A charcteristic function (CF), of a probability distribution, $\psi(t)$ is randomly self-decomposable (RSD) if for each $p,c\in [0,1]$ there exists a CF $\psi_{c,p}(t)$ such that 
\begin{equation}
\psi(t)= \psi_{c,p}(t)\{p+(1-p)\psi(ct)\}.
\end{equation}

Kozubowski and Podg\'orski (2010) then discusses the relation of RSD laws to SD laws and geometrically infinitely divisible (GID) laws. In proposition 2.3 they prove, in an elegent manner, that the class of RSD laws equals the intersection of the classes of GID laws and SD laws. They also discuss a variety of examples. We need the following in our discussion.

\vspace{2mm}

\noindent \textbf {Definition 1.2} Harris$(a,k)$ distribution on $\{1, 1+k, 1+2k, ....\}$ is described by its probability generating function (PGF)
\begin {equation}
P(s)= \frac{s}{\{a-(a-1)s^k\}^{1/k}}, \text {$k>0$ integer and $a>1$.}
\end {equation}

\vspace{2mm}

\noindent \textbf {Definition 1.3} A CF $\psi(t)$ is Harris-ID (HID) if for each $p\in(0,1)$ there exists a CF $\psi_p(t)$ such that
\begin {equation}
\psi(t)=\frac{\psi_p(t)}{\{a-(a-1)\psi_{p}^{k}(t)\}^{1/k}}, p=\frac{1}{a}.
\end {equation}
 
\vspace{2mm}

\noindent \textbf {Theorem 1.1} (Satheesh \textit{et al.} (2008)) A CF $\psi(t)$ is HID \textit{iff}
\begin{equation}
\psi(t)=\frac{1}{(1-\log h(t))^{1/k}}
\end{equation}
where $k>0$ integer and $h(t)$ is some CF that is ID. 

\vspace{2mm}

When $k=1$ Harris distribution becomes the geometric(\textit{p}) distribution on $\{1, 2, ...\}$ with $p=\frac{1}{a}$. For more on this distribution see Sandhya \textit{et al.} (2008). Certain aspects of HID laws and generalized AR(1) models have been discussed in Satheesh \textit{et al.} (2008). In section 2, the notion of RSD is generalized, its relation to SD laws and HID laws are presented and its connection to a stationary generalized AR(1) model is given. The notion is then extended to $\mathbf{Z_+}$-valued distributions in section 3. We closely follow the development in Kozubowski and Podg\'orski (2010).

\section{\large Generalizing RSD distributions}\label{sec2}

\noindent \textbf {Remark 2.1} In the paragraph after their Proposition 3.1 Kozubowski and Podg\'orski (2010) state that AR(1) processes described by (1.2) cannot be constructed with either (general) gamma or Gaussian distributions for $X_n$ as neither of them are GID although both are SD. However, it should be noted that gamma$(\alpha,\lambda)$ distributions (\textit {equation} (2.10)) are GID if $\alpha \leq 1$, \textit {see e.g.} Yannaros (1988) or Sandhya (1991).

\vspace{2mm}

\noindent \textbf {Definition 2.1} A CF $\psi(t)$ is Harris-RSD (HRSD) if for each $c\in (0,1]$ and each $p\in[0,1) $ there exists a distribution with CF $\psi_{c,p}(t)$ such that
\begin{equation}
\psi(t)= \psi_{c,p}(t)\{p+(1-p)\psi^{k}(ct)\}^{1/k}.
\end{equation}

\vspace{2mm}

\noindent \textbf {Remark 2.2} With the above nomenclature the RSD defined by Kozubowski and Podg\'orski (2010) is geometric RSD (GRSD) because it bridges the notions of SD and GID where as our definition bridges the notions of SD and HID.

\vspace{2mm}

When $p=0$ equation (2.1) reduces to
\begin{equation}
\psi(t)=\psi(ct) \psi_{c}(t)
\end{equation}
where $\psi_{c}(t)= \psi_{c,0}(t)$, that is $\psi(t)$ is SD. 
On the other hand when $c=1$ equation (2.1) becomes
\begin{equation}
\psi(t)= \psi_{p}(t)\{p+(1-p)\psi^{k}(t)\}^{1/k}.
\end{equation}
where $\psi_{p}(t)=\psi_{1,p}(t)$. Solving for $\psi(t)$ we get
\begin{equation}
\psi(t)= \frac{\psi_p(t)}{\{a-(a-1)\psi_p^k(t)\}^{1/k}} ;
a=\frac1p.
\end{equation}
That is $\psi(t)$ is HID.

\vspace{2mm}

Denoting the classes of HRSD, SD and HID distributions by $\mathcal{C}_{HRSD}$, $\mathcal{C}_{SD}$ and $\mathcal{C}_{HID}$ the above discussion shows that $\mathcal{C}_{HRSD}\subset \mathcal{C}_{SD} \cap \mathcal{C}_{HID}$. In the next Proposition we show that we have equality here.

\vspace{2mm}
\noindent \textbf {Proposition 2.1} We have $\mathcal{C}_{HRSD} = \mathcal{C}_{SD} \cap \mathcal{C}_{HID}$. Further, whenever the CF $\psi(t) \in \mathcal{C}_{HRSD}$, the CF $\psi_{c,p}(t)$ in (2.1) can be written as
\begin {equation}
\psi_{c,p}(t)=\psi_{c}(t). \psi_{p}(ct)
\end{equation}
where $\psi_{c}(t)$ and $\psi_{p}(t)$ are given by
\begin {equation}
\psi_{c}(t)= \frac {\psi(t)}{\psi(ct)}
\end{equation}
\begin {equation}
\psi_{p}(t)= \frac {\psi(t)} {\{p+(1-p)\psi^{k}(t)\}^{1/k}}   
\end{equation}

\vspace{2mm}

\noindent \textit {Proof.} If the CF $\psi(t)$ is SD then for each $c\in(0,1]$ the function $\psi_{c}(t)$ in (2.6) is a genuine CF and similarly if $\psi(t)$ is HID then for each $p\in[0,1)$ the function $\psi_{p}(t)$ in (2.7) also is a genuine CF. Consequently (2.5) is a well defined CF and hence (2.1) holds, proving the assertion.  

\vspace{2mm}

Now let us consider a generalization of the AR(1) sequence (1.2). Here $\{X_{n}\}$ is composed of $k$ independent AR(1) sequences $\{Y_{n,i}\},i=1,2, \dots\ k$ and where for each \textit{n}, $\{Y_{n,i}\}$ are independent. That is, for each $n$, $X_{n}=\sum_{i=1}^k Y_{n,i}$ and $\epsilon_{n}=\sum_{i=1}^k \epsilon_{n,i}$ where $\{Y_{n,i}\}$ is an \textit{i.i.d} sequence and similarly $\{\epsilon_{n,i}\}$ is also an \textit{i.i.d} sequence, $k$ being a fixed positive integer. Further, it is also assumed that for each $n$, $\epsilon_{n,i}$ is independent of $Y_{n-1,i}$ for all $i=1,2, \dots\ k$. Situations where such a model can be useful have been discussed in Satheesh \textit{et al.} (2008). 
\begin{equation}
\sum_{i=1}^k Y_{n,i}=
\begin{cases}
\sum_{i=1}^k \epsilon_{n,i}, \text{ with probability }p,\\
\sum_{i=1}^k cY_{n-1,i} + \sum_{i=1}^k \epsilon_{n,i},  \text{ with probability}(1-p).
\end{cases}
\end{equation}

\vspace{2mm}

In terms of CFs and assuming stationarity we get
\begin{equation}
\psi_Y(t) = \psi_\epsilon (t) \{p+(1-p)\psi_Y^k (ct)\}^{1/k}.
\end{equation}

The following Proposition is now clear.

\vspace{2mm}

\noindent \textbf {Proposition 2.2\textit{a}} If $\{Y_{n,i}\}$ describes the AR(1) model (2.8) that is stationary for each $c\in (0,1]$ and $p\in[0,1)$ then the distribution of $\{Y_{n,i}\}$ is HRSD. 

Now proceeding as in the proof of Theorem 2.2 in Bouzar and Satheesh (2008) we have the following converse to Proposition 2.2\textit{a}.

\vspace{2mm}

\noindent \textbf {Proposition 2.2\textit{b}} If $\psi_{Y}(t)$ is a CF that is HRSD with $\psi_{c,p}(t)=\psi_\epsilon(t)$ for each $c\in (0,1]$ and $p\in[0,1)$ then there exists a stationary AR(1) model described by (2.8) with $\psi_{Y}(t)$ the CF of $\{Y_{n,i}\}$ and $\psi_\epsilon(t)$ that of the innovations $\{\epsilon_{n,i}\}$.

\vspace{2mm}

\noindent \textbf {Example 2.1}
Gamma$(\frac{1}{k},\lambda)$ distributions has CF
\begin {equation}
\psi(t)=\frac{1}{(1-i \lambda t)^{1/k}}, \text {$k>0$ integer, $\lambda >0$}
\end{equation}
is HID. Further since it is also SD, this distribution is HRSD. 

\vspace{2mm}

\noindent \textbf {Example 2.2} Let the CF $\psi(t)$ be Harris-sum-stable for every $c\in(0,1)$. Then
\begin {equation}
\psi(t)=\psi(ct).\frac{1}{\{a-(a-1)\psi^k(ct)\}^{1/k}}.
\end{equation}

\noindent The second factor on the RHS is also a genuine CF being a Harris-sum of $\psi(t)$ where the support of this Harris distribution is $\mathbf{Z_+}$. Thus $\psi(t)$ is SD. Hence if we take $h(t)$ as a stable CF in Theorem 1.1 then $\psi(t)$ in (1.6) is Harris-sum-stable for every $c\in(0,1)$ and we have a general procedure to construct CFs that are HRSD. With $k=1$ above, we have the corresponding geometric-sum-stable laws and a procedure to construct the examples in Kozubowski and Podg\'orski (2010).

\section{\large Discrete analogue of HRSD distributions}\label{sec3}
Steutel and van Harn (1979) had developed discrete SD (DSD) distributions. We now introduce RSD and HRSD for $\mathbf{Z_+}$-valued distributions. Some aspects of discrete HID (DHID) laws and generalized AR(1) models on $\mathbf{Z_+}$ have been discussed in Satheesh \textit{et al.} (2010\textit{b}). 
\vspace{4mm}

\noindent \textbf {Definition 3.1} A PGF $P(s)$ is DHID if for each $p\in(0,1)$ there exists a PGF $P_p(s)$ such that
\begin {equation}
P(s)=\frac{P_p(s)}{\{a-(a-1)P_{p}^{k}(s)\}^{1/k}}, p=\frac{1}{a}.
\end {equation}

\vspace{2mm}

\noindent \textbf {Theorem 3.1} (Satheesh \textit{et al.} (2010\textit{a})) A PGF $P(s)$ is DHID \textit{iff} 
\begin {equation}
P(s)=\frac{1}{(1-\log R(s))^{1/k}},
\end {equation}
where $k>0$ integer and $R(s)$ is a PGF that is DID.

\vspace{4mm}

\noindent \textbf {Definition 3.2} A PGF $P(s)$ is discrete HRSD (DHRSD) if for each $c\in (0,1]$ and $p\in[0,1) $ there exists a PGF $P_{c,p}(s)$ such that
\begin{equation}
P(s)= P_{c,p}(s)\{p+(1-p)P^{k}(1-c+cs)\}^{1/k}.
\end{equation}

\vspace{4mm}
 
Denoting the classes of DHRSD, DSD and DHID distributions by $\mathcal{C}_{DHRSD}$, $\mathcal{C}_{DSD}$ and $\mathcal{C}_{DHID}$ we can proceed as in Section 2 to arrive at

\vspace{2mm}

\noindent \textbf {Proposition 3.1} We have $\mathcal{C}_{DHRSD} = \mathcal{C}_{DSD} \cap \mathcal{C}_{DHID}$. Further, whenever $P(s) \in \mathcal{C}_{DHRSD}$, the PGF $P_{c,p}(s)$ in (3.3) can be written as
\begin {equation}
P_{c,p}(s)=P_{c}(s). P_{p}(1-c+cs)
\end{equation}
where $P_{p}(s)$ and $P_{c}(s)$ are given by
\begin {equation}
P_{c}(s)= \frac {P(s)}{P(1-c+cs)}
\end{equation}
\begin {equation}
P_{p}(s)= \frac {P(s)} {\{p+(1-p)P^{k}(s)\}^{1/k}}   
\end{equation}

\vspace{2mm}

Again, considering the $\mathbf{Z_+}$-valued analogue of the generalized AR(1) scheme (2.8) with $\odot$, the binomial thinning operator in Steutel and van Harn (1979) we have the INAR(1) model
\begin{equation}
\sum_{i=1}^k Y_{n,i}=
\begin{cases}
\sum_{i=1}^k \epsilon_{n,i}, \text{ with probability }p,\\
\sum_{i=1}^k c \odot Y_{n-1,i} + \sum_{i=1}^k \epsilon_{n,i},  \text{ with probability}(1-p).
\end{cases}
\end{equation}

\vspace{2mm}

Assuming stationarity we have the following Propositions as in Section 2.
\vspace{2mm}

\noindent \textbf {Proposition 3.2\textit{a}} If $\{Y_{n,i}\}$ describes the stationary INAR(1) model (3.7) for each $c\in (0,1]$ and $p\in[0,1)$ then the distribution of $\{Y_{n,i}\}$ is DHRSD.
\vspace{2mm}

Conversely,

\vspace{2mm}

\noindent \textbf {Proposition 3.3\textit{b}} If $P_{Y}(s)$ is a PGF that is DHRSD with $P_{c,p}(s)= P_\epsilon(s)$ for each $c\in (0,1]$ and $p\in[0,1)$ then there exists a stationary INAR(1) model described by (3.7) with $P_{Y}(s)$ the PGF of $\{Y_{n,i}\}$ and $P_\epsilon(s)$ that of $\{\epsilon_{n,i}\}$.

\vspace{2mm}

\noindent \textbf {Example 3.1} Negative binomial$(\frac{1}{k},\lambda)$ distributions with PGF
\begin {equation}
P(s)=\frac{1}{(1+\lambda(1-s))^{1/k}}, \text {$k>0$ integer, $\lambda >0$} 
\end {equation} 
is DHID. Further since it is also DSD, this distribution is HRSD.

\vspace{2mm}
\noindent \textbf {Example 3.2} We may also proceed in a general frame work as done in Example 2.2 to construct PGFs that are HRSD.

\vspace{2mm}

Satheesh and Sandhya (2010) has proposed a further generalization of HRSD distributions based on the notion of $\mathcal{N}ID$ distributions of Gnedenko and Korolev (1996). 
%%%%%%%%%%%%%%%%%%%%%%%%%%%%%%%%%%%%%%%%%%%

\end{document}